\documentclass[a4paper,12pt]{amsart}

\usepackage{amsfonts,amssymb,amsthm,amscd}

\newtheorem{teo}{Theorem}[section]
\newtheorem{lema}[teo]{Lemma}
\newtheorem{prop}[teo]{Proposition}
\newtheorem{cor}[teo]{Corollary}
\newtheorem{obs}[teo]{Remark}
\newtheorem{subsct}[teo]{}
\theoremstyle{plain}

\newcommand{\bdem}{\begin{proof}}
\newcommand{\edem}{\end{proof}}

\newcommand{\pr}{\mathbb{P}}
\newcommand{\inteiro}{\mathbb{Z}}
\newcommand{\complexo}{\mathbb{C}}
\newcommand{\reg}{\mathrm{reg}}
\renewcommand{\O}{\mathcal O}
\newcommand{\I}{\mathcal I}
\newcommand{\F}{\mathcal F}
\renewcommand{\L}{\mathcal L}

\renewcommand{\S}{\mathcal S}
\newcommand{\w}{\omega}
\newcommand{\Imag}{\mathrm{Im}}
\newcommand{\lra}{\longrightarrow}
\renewcommand{\:}{\colon}
\newcommand{\wt}{\widetilde}
\newcommand{\ox}{\otimes}
\newcommand{\T}{\mathcal T}

\begin{document}
\author[Cruz and Esteves]
{Joana D. A. S. Cruz \and Eduardo Esteves}
\thanks{First author supported by CAPES and FAPEMIG, Processo 428/07, second author supported by CNPq, Processos 300004/95-8 and 470761/2006-7.}
\title[Regularity of subschemes invariant under Pfaff fields on $\pr^n_k$]
{Bounding the regularity of subschemes invariant under Pfaff fields on projective spaces}
\begin{abstract} 
A Pfaff field on $\pr^n_k$ is a map $\eta\:\Omega^s_{\pr^n_k}\to\L$ from the sheaf of 
differential $s$-forms to an invertible sheaf. The interesting ones are those 
arising from a Pfaff system, as they give rise to a distribution away from their singular locus. 
A subscheme $X\subseteq\pr^n_k$ is said to be invariant under $\eta$, if $\eta$ induces a Pfaff 
field $\Omega^s_X\to\L|_X$. We give bounds for the Castelnuovo--Mumford regularity of invariant 
complete intersection subschemes (more generally, arithmetically Cohen--Macaulay subschemes) of dimension $s$, depending 
on how singular these schemes are, thus bounding the degrees of the hypersurfaces that cut them out. 
\end{abstract}
\maketitle

\section{Introduction}

In 1891, Poincar\'e \cite{Po}, p.~161, posed the problem of bounding a priori the 
degree of the first integral of a polynomial vector field on the complex plane, 
when the integral is algebraic. The importance of such a bound is that 
it allows us to decide whether the integral is algebraic or not by 
making purely algebraic computations.

Poincar\'e himself produced bounds in special cases. But no bounds have 
been found in general. Actually, many obstructions to finding such bounds 
have been discovered: For instance, Lins Neto \cite{L} produced examples 
to show that a bound cannot depend only on the degree $m$ of the 
vector field and on the analytic type of its singularities in the plane or at infinity.

The current interest in Poincar\'e's problem was revived exactly a 
hundred years later by Lins Neto and Cerveau \cite{CeL}, who showed 
that 
an algebraic curve invariant under the vector field has degree at most 
$m+2$, if the singularities of the curve are ordinary double points, the 
bound achieved only if the curve is reducible; see loc.~cit., Thm.~1, p.~891. 
Since then many papers 
have concentrated on this related problem, of bounding 
the degrees of algebraic curves invariant under the vector field. This 
has often been called the \emph{Poincar\'e problem}. Works on this 
problem, allowing for more singular curves, are \cite{CamCar}, 
\cite{Car}, \cite{EsK3}, \cite{dPW} and \cite{Pe}, to cite a few. 

The problem has also been considered for higher dimensional spaces. One of the 
first to do so was Soares \cite{S}. 
In today's language, and in great generality, let $\pr^n_k$ denote the 
$n$-dimensional projective space over an algebraically closed field $k$, 
and consider a \emph{Pfaff field}, a map 
$\eta\:\Omega^s_{\pr^n_k}\to\L$ from the sheaf of $s$-forms 
$\Omega^s_{\pr^n_k}:=\bigwedge^s\Omega^1_{\pr^n_k}$, for an integer 
$s$ between 1 and $n-1$ called the \emph{rank} of $\eta$, 
to an invertible sheaf $\L$. The unique 
numerical global invariant under deformations 
of $\eta$ besides its rank is $m:=\deg(\L)+s$, the 
\emph{degree} of $\eta$. 
The \emph{singular locus} 
of $\eta$ is its degeneracy scheme $\S$, supported on the set 
of points where $\eta$ is not surjective. 
A closed subscheme $X\subseteq\pr^n_k$ is said to be 
\emph{invariant} under $\eta$ if $\eta$ induces a Pfaff field 
$\Omega^s_X\to\L|_X$ on $X$. The above terminology is taken from \cite{EsK2}, Section 3.

The Pfaff field $\eta$ may arise from a Pfaff system, as defined by 
Jouanolou \cite{J}, pp.~136--138, which is essentially a ``singular distribution,'' in particular 
giving rise to an actual distribution on $\pr^n_k-\S$. 
Then subschemes of pure dimension $s$ that are invariant under $\eta$ are 
solutions of the corresponding Pfaff systems; see \cite{EsK2}, Prop.~3.2, p.~3782 for a precise 
statement.

If $s=1$ then $\eta$ is the homogenization of a polynomial 
vector field on $\mathbb C^n$. If $s=n-1$, through the perfect pairing 
$\Omega^s_{\pr^n_k}\ox\Omega^{n-s}_{\pr^n_k}\to\Omega^n_{\pr^n_k}$, we 
may view $\eta$ as the homogenization of a polynomial differential 1-form 
on $\mathbb C^n$. In both cases, $\eta$ arises from a distribution away 
from $\S$.

Some of the statements in the literature, and all of the statements in 
the present article, work in positive characteristic, under suitable 
assumptions. However, to simplify the ongoing discussion, assume 
that $k$ has characteristic zero.

For $s=n-1$ one may search for bounds on the degrees of hypersurfaces 
invariant under $\eta$. For instance, 
under the harmless assumption that $\dim(\S)\leq n-2$, Brunella and Mendes \cite{BMe}  
showed that an invariant reduced hypersurface with at most 
normal-crossings singularities has degree at most $m+2$, 
generalizing the theorem by Cerveau and Lins Neto mentioned above; see 
loc.~cit., p.~594 for a more general statement.

For $s=1$ many inequalities have been produced for the degree and the 
genus of (reduced, equidimensional) curves invariant under 
$\eta$, for instance in \cite{CamCarG} and \cite{EsK1}. However, in the 
spirit of Poincar\'e's original problem, one should look for bounds on 
global invariants that could reduce to purely algebraic computations the 
question of whether $\eta$ has an invariant curve or not. 
The (Castelnuovo--Mumford) regularity is such an invariant, as it is well-known 
that a subscheme $X\subseteq\pr^n_k$ is cut out by hypersurfaces with degree at most its 
regularity, $\reg(X)$.

In \cite{Es} the second author shows that 
an invariant arithmetically Cohen--Macaulay (a.C.M.) 
curve $C$, with at most ordinary 
double points for singularities, such that $\S\cap C$ is finite 
has regularity at most $m+2$, with 
equality only if the curve is reducible; see loc.~cit., Thm.~1, p.~3. 
Since complete intersections are 
a.C.M., and since the regularity of a 
hypersurface is its degree, the statement is another 
generalization of Cerveau's and Lins Neto's result.

Later, the second author and Kleiman showed that the 
inequality $\reg(X)\leq m+2$ for an invariant a.C.M.~curve 
(for $s=1$) or invariant reduced hypersurface $X\subseteq\pr^n_k$ (for $s=n-1$) 
with normal-crossings singularities was a consequence of the fact that 
$h^s(\Omega^s_X(1))=0$, and that the same holds for intermediate $s$. 
More precisely, for any $s$, an invariant, reduced, 
a.C.M.~subscheme $X\subseteq\pr^n_k$ 
of pure dimension $s$ whose irreducible components are not contained 
in $\S$ has regularity bounded by $m+2$ if $h^s(\Omega^s_X(1))=0$, and 
bounded by $m+1$ if $h^s(\Omega^s_X)=1$; see \cite{EsK2}, Cor.~4.5, p.~3790 and Rmks.~4.6 and 
4.7, p.~3791, from which the assertion can be extracted. 

However, no further conditions for 
when $h^1(\Omega^1_X)=1$ or $h^1(\Omega^1_X(1))=0$ are given in \cite{EsK2}. 
These appear later in \cite{EsK3}, by the same authors, but only for $n=2$. There a (reduced) 
plane curve $C$ of degree $d$ is considered, and it is shown that if the singular locus of $C$ 
has regularity $\sigma$ bounded by 
$d-2$ then $h^1(\Omega^1_C)=1$; and hence $d\leq m+1$ if 
$C$ is invariant. The highly singular case is handled as well, being shown that if $C$ is 
invariant and $\rho:=\sigma-d+2$ is 
positive, then $d\leq m+1+\rho$, with equality if $d\geq 2m+2$ and $\S$ 
is finite; see loc.~cit., Thm.~2.5, p.~61.

In the present article, we extend the results of \cite{EsK3} for $n>2$ 
and any $s$. More precisely, our Theorem \ref{h11} states that a connected, reduced 
subscheme $X\subseteq\pr^n$ 
of pure dimension $s>0$ satisfies $h^s(\Omega^s_X)=1$ if 
$X$ is a.C.M. and subcanonical, for instance a complete intersection, 
and if its singular locus has regularity $\sigma$ bounded by $r-2$, where 
$r$ is the regularity of $X$. From it follows Theorem \ref{principal0}, stating that 
$r\leq m+1$ if in addition 
$X$ is invariant and $\dim(\S\cap X)<s$. Furthermore, by our Theorem \ref{principal1}, 
if $X$ is simply a.C.M., and is invariant with 
$\dim(\S\cap X)<s$, then $r\leq m+1+\rho$, where 
$\rho:=\max(1,\sigma-r+2)$. Finally, Theorem \ref{igualdade1} says that $r=m+1+\rho$ if all the 
following conditions hold: $s=1$ and $\S$ is finite; $X$ is a.C.M., 
subcanonical and invariant; $r\geq 5$ if $m=1$ or $r\geq mn-n+4$ if $m>1$. 

Since complete intersections are a.C.M., subcanonical subschemes, we obtain as a corollary that, if 
$X\subseteq\pr^n_k$ is a reduced complete intersection of hypersurfaces of degrees $d_1,\dots,d_{n-s}$, 
and is invariant under $\eta$ with $\dim(\S\cap X)<s$, then
$$
d_1+\dots+d_{n-s}\leq\left\{\begin{array}{lc}
 m+n-s,       & \mbox{ if } \rho\leq 0,  \\
 m+n-s+\rho,  & \mbox{ if } \rho > 0,
\end{array}\right.
$$
where $\rho:=\sigma+n-s+1-d_1-\cdots-d_{n-s}$, with $\sigma$ denoting the regularity of the singular 
locus of $X$; see Corollary \ref{regfolcur}.

The techniques we use are quite simple: basically, a detailed analysis of the long exact sequences in 
cohomology of several short exact sequences of sheaves associated to the problem. 

The pervasive hypothesis of arithmetic Cohen--Macaulayness is necessary, as the example of a sequence of 
smooth curves in $\pr^3$ of increasing regularity but invariant under degree-1 rank-1 Pfaff fields, 
presented in \cite{Es}, Rmk.~21, p.~14, shows. What is not clearly necessary is the hypothesis of 
subcanonicalness.

The possibility that $r=m+1+\rho$ is investigated only for $s=1$, because then $\S$ is easier to 
understand. Then, if $\S$ has dimension 0, which is the expected dimension and the case when $\eta$ is 
general, the regularity of $\S$ is 
1 if $m=1$ and $mn-n+2$ if $m>1$; see Proposition \ref{regsing} and the remark thereafter. This 
regularity gives the bound above which $r$ must be for the equality $r=m+1+\rho$ to hold. 
On the other hand, for $s\geq 2$, those $\eta$ having an invariant reduced subscheme of 
pure dimension $s$ have large singular locus; indeed, $\dim(\S)\geq s-1$ by 
\cite{EsK2}, Cor.~4.5, p.~3790. In particular, $\S$ does no have the expected dimension.

Section 2 collects a few results on the Castelnuovo--Mumford regularity and on 
arithmetically Cohen--Macaulay subschemes. In Section 3 we give conditions for when a subscheme 
$X\subseteq\pr^n_k$ of pure dimension $s$ satisfies $h^s(\Omega^s_X)=1$. In Section 4 we prove our 
bounds on the regularity of closed subschemes invariant under Pfaff fields. Finally, in Section 5 
we prove that these bounds are attained, if the regularity is large enough, in the case of rank-1 Pfaff 
fields.

\section{Arithmetically Cohen--Macaulay subschemes} 

\begin{subsct}\normalfont (\textit{The Castelnuovo--Mumford regularity}) 
Fix a positive integer $n$. Given $m\in\inteiro$, 
we say that a coherent sheaf $\F$ on $\pr^n_k$ is \textit{$m$-regular} if
$H^i(\F(m-i))=0$ for each integer $i>0$. 

Let $X\subseteq\pr^n_k$ be a closed subscheme. If $X\neq\pr^n_k$ then the
\textit{Castelnuovo--Mumford regularity} of $X$, or simply \textit{regularity}, 
is the smallest integer $m$ for which its sheaf of ideals 
is $m$-regular. By definition, the regularity of $\pr^n_k$ is 
$1$. Denote the regularity of $X$ by $\reg(X)$. 

The regularity is well-defined. In fact, let $\I_X$ denote the sheaf of ideals 
of $X$, and consider the natural exact sequence:
\begin{equation}\label{XPN}
0\to\I_X\to\O_{\pr^n_k}\to\O_X\to 0.
\end{equation}
Twisting it by $m-n$ and taking cohomology we get the following exact sequence:
$$
H^n(\I_X(m-n))\lra H^n(\O_{\pr^n_k}(m-n))\lra H^n(\O_X(m-n)).
$$
The middle group is zero if and only if $m\geq 0$. If $X\neq\pr^n_k$ then the last 
group is zero, and hence $H^n(\I_X(m-n))=0$ only if $m\geq 0$. 

The above reasoning shows that $\reg(X)\geq 0$. Furthermore, $\reg(X)=0$ if 
and only if $X=\emptyset$. Indeed, if $X$ is empty, $\I_X=\O_{\pr^n_k}$, which is 
$0$-regular by Serre computation. On the other hand, if $\I_X$ is 0-regular then 
$\I_X$ is globally generated, by \cite{Mu}, p.~99. Since $\reg(X)\neq 1$, we have 
that $X\neq\pr^n_k$, and hence $\I_X\neq 0$. So $H^0(\I_X)\neq 0$, which 
implies that $\I_X=\O_{\pr^n_k}$, and thus $X=\emptyset$.

Also, $\reg(X)=1$ if and only if $X$ is a linear subspace of $\pr^n_k$. Indeed, 
if $\reg(X)=1$ then $\I_X(1)$ is globally generated, which implies that $X$ is 
cut out by a system of hyperplanes. Conversely, suppose $X$ is a linear subspace 
of $\pr^n_k$. Twisting (\ref{XPN}) by $1-i$ and taking cohomology, we get the 
following exact sequence:
$$
H^{i-1}(\O_{\pr^n_k}(1-i))\lra H^{i-1}(\O_X(1-i))\lra H^i(\I_X(1-i))\lra 
H^i(\O_{\pr^n_k}(1-i)).
$$
If $i>1$, the second and last groups are zero, by Serre computation. Thus 
$H^i(\I_X(1-i))=0$ for $i>1$. For $i=1$ the last group is zero, and the first map 
is an isomorphism. Thus $H^1(\I_X)=0$. So $\reg(X)\leq 1$. Since $X\neq\emptyset$, 
it follows that $\reg(X)=1$.
\end{subsct}

\begin{prop}\label{regconjfin} Let $X\subseteq\pr^n_k$ be a
closed subscheme. If $\dim(X)=0$ then 
$\reg(X)$ is the smallest nonnegative integer $r$ such that $H^1(\I_X(r-1))=0$, 
where $\I_X$ is the sheaf of ideals of $X$.
\end{prop}

\bdem Clearly, $H^i(\I_X(m))=0$ for every $i>n$ and every $m\in\inteiro$. 
Thus the assertion follows from the definition of 
regularity if $n=1$. 

Suppose now that $n>1$. We need only show that $H^i(\I_X(r-i))=0$ for each $i=2,\dots,n$ and each 
$r\geq 0$. Let $m\in\inteiro$. Since $X$ has dimension zero, $H^i(\O_X(m))=0$ for every 
$i\geq 1$. 
On the other hand, from Serre computation, $H^i(\O_{\pr^n_k}(m))=0$
for each $i=1,\ldots, n-1$. Twisting the natural
exact sequence
\begin{equation}\label{seid}
0\rightarrow \I_X\rightarrow \O_{\pr^n_k}\rightarrow \O_X\rightarrow 0
\end{equation}
by $m$, and taking cohomology, we get, for each $i=2,\ldots, n$, the exact sequence
$$
H^{i-1}(\O_X(m))\rightarrow H^i(\I_X(m))\rightarrow H^i(\O_{\pr^n_k}(m))\rightarrow 
H^i(\O_X(m)).
$$ 
If $i=2,\ldots,n-1$ then $H^{i-1}(\O_X(m))=H^i(\O_{\pr^n_k}(m))=0$, 
and hence $H^i(\I_X(m))=0$. If $i=n$, since $H^{n-1}(\O_X(m))=H^n(\O_X(m))=0$ because $n\geq 2$, we 
have
$$
H^n(\I_X(m))\cong H^n(\O_{\pr^n_k}(m)).
$$
But, from Serre computation, $H^n(\O_{\pr^n_k}(m))=0$ if $m\geq -n$. Thus $H^n(\I_X(r-n))=0$ for 
each $r\geq 0$.
\edem

\begin{subsct}\normalfont (\textit{Arithmetically Cohen--Macaulay subschemes}) An 
equidimensional closed subscheme $X\subseteq\pr^n_k$ is said to be 
\textit{arithmetically Cohen--Macaulay} (or simply {\rm a.C.M.}) if its 
coordinate ring is Cohen--Macaulay. Alternatively, if $X$ has positive dimension, 
$X$ is a.C.M. if the restriction map
$$
H^0(\O_{\pr^n_k}(m))\lra H^0(\O_X(m))
$$
is surjective and $H^j(\O_X(m))=0$ for each $m\in\inteiro$ and 
$j=1,\dots,\dim(X)-1$. Or,equivalently, $X$ is a.C.M. if 
$H^j(\I_X(m))=0$ for each $m\in\inteiro$ and $j=1,\dots,\dim(X)$,
where $\I_X$ is the sheaf of ideals of $X$. Notice that it follows that 
$h^0(\O_X)=1$, and hence that $X$ is connected.

Complete intersections are the simplest examples of a.C.M.~subschemes.
\end{subsct}

\begin{prop}\label{reginter1}
Let $X\subseteq\pr^n_k$ be a closed subscheme of pure
dimension $s>0$. If $X$ is arithmetically Cohen--Macaulay then $\reg(X)$ 
is the smallest nonnegative integer $r$ such that 
$H^s(\O_X(r-s-1))=0$.
\end{prop}

\bdem  Suppose first that $s=n$, that is, $X=\pr^n_k$. By
definition, the regularity of $\pr^n_k$ is 1. On the other hand, by 
Serre computation, $H^n(\O_{\pr^n_k}(r-n-1))=0$ if and
only if $r\geq 1$. So, the proposition holds for $s=n$.

Now, assume $s<n$. Let $\I_X$ denote the sheaf of ideals of $X$. 
Since $X$ is a.C.M., $H^i(\I_X(r-i))=0$ for every $r\in\inteiro$ and 
each $i=1,\ldots, s$. 
On the other hand, twisting the natural short exact sequence 
\begin{equation}\label{sqeideal} 
0\to\I_X\to\O_{\pr^n_k}\to\O_X\to0
\end{equation}
by $r-i$, and taking cohomology, we get the following 
exact sequence, for each integer $i>0$:
\begin{equation}\label{sqi}
H^{i-1}(\O_{\pr^n_k}(r-i))\lra H^{i-1}(\O_X(r-i))\lra H^i(\I_X(r-i))\lra 
H^i(\O_{\pr^n_k}(r-i)).
\end{equation}
For $i=s+2,\ldots,n-1$, since 
$$
H^{i-1}(\O_X(r-i))= H^i(\O_{\pr^n_k}(r-i))=0,
$$ 
we have that $H^i(\I_X(r-i))=0$ for every $r\in\inteiro$. Also, since 
$H^n(\O_{\pr^n_k}(r-n))=0$ for $r\geq 0$, and $H^{n-1}(\O_X(r-n))=0$ if $s<n-1$, 
it follows that $H^n(\I_X(r-n))=0$ for every $r\geq 0$, if $s<n-1$.

So, $\reg(X)$ is the smallest nonnegative integer $r$ such that 
$H^{s+1}(\I_X(r-s-1))=0$. But, if $r\geq 0$ then
$$
H^s(\O_{\pr^n_k}(r-s-1))=H^{s+1}(\O_{\pr^n_k}(r-s-1))=0,
$$
because $0<s<n$, by Serre computation. So, by the exactness of (\ref{sqi}) for $i=s+1$,
$$
H^s(\O_X(r-s-1))\cong H^{s+1}(\I_X(r-s-1))
$$ 
for every integer $r\geq 0$.
\edem

\begin{subsct}\normalfont (\textit{Subcanonical subschemes}) 
Let $X\subseteq\pr^n_k$ be a closed subscheme. Let $\w_X$ be the 
dualizing sheaf of $X$, that is, 
$$
\w_X:=\mathcal{E}xt^{n-s}_{\O_{\pr^n_k}}(\O_X,\O_{\pr^n_k}(-1-n)),
$$
where $s:=\dim(X)$. If there is $a\in \inteiro$
such that $\omega_X\cong\O_X(a)$, then we say that $X$ is
$a$-\textit{subcanonical} (or simply \textit{subcanonical}). If $\dim(X)>0$ 
then $a$ is unique.
\end{subsct}

\begin{prop}\label{regsub}
Let $X\subseteq \pr^n_k$ be an arithmetically Cohen--Macaulay 
$a$-subcanonical subscheme of 
pure dimension $s>0$. Then $a\geq -s-1$ and $\reg(X)=a+s+2$.
\end{prop}

\bdem Observe first that $H^0(\O_X(i))=0$ if and only if
$i<0$. Indeed, since $\O_X(1)$ is very ample, $h^0(\O_X(i))>0$ for $i\geq 0$. On 
the other hand, if there were a nonzero global section of $\O_X(i)$ 
for a certain $i<0$, 
multiplying it by $-iH$, for a sufficiently general hyperplane 
section $H\subset X$, we 
would obtain a nonzero global section of $\O_X$ vanishing at $H$, an absurd.
 
By duality, since $\O_X(a)$ is the dualizing sheaf of $X$,
$$
h^s(\O_X(r-s-1))=h^0(\O_X(a-r+s+1)).
$$
So, $H^s(\O_X(r-s-1))=0$ if and only if
$a-r+s+1<0$, that is, if and only if $r\geq a+s+2$. It follows now from
Proposition~\ref{reginter1} that $\reg(X)=\max(a+s+2,0)$. However, since 
$\reg(X)>0$, we have that $a+s+2>0$ and $\reg(X)=a+s+2$.
\edem

\section{The singular locus} 

\begin{subsct}\label{canonical}\normalfont (\textit{The singular locus of a $k$-scheme}) 
Let $X$ be an algebraic $k$-scheme. For each integer $s\geq 0$, denote by 
$\Omega^s_X$ the sheaf of K\"ahler $s$-forms of $X$, that is, 
$\Omega^s_X:=\bigwedge^s\Omega^1_X$, where 
$\Omega^1_X$ is the sheaf of K\"ahler differentials of $X$. 

Assume $X$ is reduced, projective and of pure dimension $s>0$. Let 
$\omega_X$ be its dualizing sheaf and $\gamma_X\:\Omega^s_X\to\omega_X$ the 
canonical map. The map $\gamma_X$ is constructed as follows. Let $X_1,\dots,X_m$ 
be the irreducible components of $X$ with their reduced induced subscheme structures. For 
each $i=1,\dots,m$ there is a natural map $\gamma_i\:\Omega^s_{X_i}\to\wt\w_{X_i}$, 
where $\wt\w_{X_i}$ is Kunz's sheaf of regular differential forms of $X_i$. Also, the map 
is an isomorphism on the smooth locus of $X_i$; see \cite{Ku}, pp.~103--105. 
Furthermore, by \cite{Ku}, Satz 2.2, p.~95 or \cite{Lp}, Thm.~0.2B, p.~15, the 
sheaf $\wt\w_{X_i}$ is dualizing, in a natural way; so there is a natural isomorphism 
$\xi_i\:\wt\w_{X_i}\to\w_{X_i}$. 

The restriction map 
$\tau\:\O_X\to\O_{X_1}\oplus\cdots\oplus\O_{X_m}$ induces 
a map $\tau'\:\w_{X_1}\oplus\cdots\oplus\w_{X_m}\to\w_X$. As $\tau$ is an isomorphism on 
the smooth locus of $X$, so is $\tau'$. Then $\gamma_X$ is, by definition, the composition
$$
\begin{CD}
\Omega^s_X @>>> \bigoplus_{i=1}^m\Omega^s_{X_i} @>(\gamma_1,\dots,\gamma_m)>> 
\bigoplus_{i=1}^m\wt\w_{X_i} @>(\xi_1,\dots,\xi_m)>> 
\bigoplus_{i=1}^m\w_{X_i} @>\tau'>> \w_X,
\end{CD}
$$
where the first map is induced by restriction. All the above maps are isomorphisms on the smooth locus 
of $X$, and thus so is $\gamma_X$.

Let $\Sigma_X$ be the scheme-theoretic support of the cokernel of 
$\gamma_X$. We call $\Sigma_X$ the \textit{singular locus} of $X$. 
Since $X$ is reduced, whence generically smooth, 
$\gamma_X$ is generically an isomorphism, and hence 
$\dim(\Sigma_X)<s$.

The sheaf $\w_X$ is torsion-free, rank-1. Indeed, it is generically isomorphic to $\Omega^s_X$, whence has 
rank 1. Its torsion subsheaf $\T(\w_X)$ is supported on a subscheme of dimension less than $s$, and 
hence $H^s(\T(\w_X))=0$. On the other hand, the injection $\T(\w_X)\to\w_X$ corresponds by duality 
to a map $H^s(\T(\w_X))\to k$. Since this map is zero, so is the injection, that is, $\T(\w_X)=0$. 

Since $\w_X$ is torsion-free, and $\gamma_X$ is generically an isomorphism, the 
kernel of $\gamma_X$ is the torsion subsheaf $\T(\Omega^s_X)\subseteq\Omega^s_X$. Thus, we get an 
injection
\begin{equation}\label{singimp}
\I_{\Sigma_X,X}\,\w_X\hookrightarrow\frac{\Omega^s_X}{\T(\Omega^s_X)}.
\end{equation}
If $X$ is Gorenstein then $\omega_X$ is invertible, and hence \eqref{singimp} is an isomorphism. 
\end{subsct}

\begin{teo}\label{h11} Let $X\subseteq \pr^n_k$ be a connected, reduced,
arithmetically Cohen--Macaulay subcanonical 
subscheme of pure dimension $s>0$. Let $\Sigma_X$ be its singular locus. 
Let $r:=\reg(X)$ and $\sigma:=\reg(\Sigma_X)$. If $\sigma=0$ or $\sigma\leq r-2$ 
then $H^s(\Omega^s_X)\cong k$.
\end{teo}

\bdem The assertion follows from Serre computation if $s=n$. Assume $s<n$. Consider the 
injection
$$
\I_{\Sigma_X,X}\,\w_X\hookrightarrow\frac{\Omega^s_X}{\T(\Omega^s_X)}.
$$
Since $X$ is reduced, both the source and target of this injection are of rank 1. So the 
injection is generically an isomorphism. Since the torsion subsheaf 
$\T(\Omega^s_X)$ is supported in dimension at 
most $s-1$, it follows that
$$
H^s(\Omega^s_X)\cong H^s(\I_{\Sigma_X,X}\,\omega_X).
$$
Since $\omega_X\cong\O_X(r-s-2)$ by Proposition \ref{regsub}, we must show that
$$
H^s(\I_{\Sigma_X,X}(r-s-2))\cong k.
$$

Set $a:=r-s-2$. Let $\I_{\Sigma_X}$ and $\I_X$ be the sheaves of ideals of $\Sigma_X$ and $X$ in 
$\pr^n_k$. We claim that 
\begin{equation}\label{s+1a}
H^{s+1}(\I_X(a))\cong k.
\end{equation}
Indeed, twisting the natural exact sequence
$$
0\to\I_X\to\O_{\pr^n_k}\to\O_X\to0
$$
by $a$ and taking cohomology, we get the exact sequence
$$
H^s(\O_{\pr^n_k}(a))\lra H^s(\O_X(a))\lra H^{s+1}(\I_X(a))\lra H^{s+1}(\O_{\pr^n_k}(a)).
$$
The first and last groups above are zero because $s<n$ and $r>0$, respectively. Thus
$$
H^{s+1}(\I_X(a))\cong H^s(\O_X(a)).
$$
But $\omega_X\cong\O_X(a)$. So, by Serre Duality,
$$
H^s(\O_X(a))\cong H^0(\O_X)\cong k,
$$
where the last isomorphism follows from the connectedness of $X$.

Now, twisting the natural exact sequence
$$
0\to\I_X\to\I_{\Sigma_X}\to\I_{\Sigma_X,X}\to0
$$
by $a$, and taking cohomology, we get the exact sequence:
\begin{equation}\label{midmap}
H^s(\I_X(a))\to H^s(\I_{\Sigma_X}(a))\to H^s(\I_{\Sigma_X,X}(a))\to H^{s+1}(\I_X(a))\to 
H^{s+1}(\I_{\Sigma_X}(a)).
\end{equation}
Since $X$ is a.C.M.~of dimension $s$, the first group is zero. The last 
group is also zero. Indeed, twisting the natural exact sequence
$$
0\to\I_{\Sigma_X}\to\O_{\pr^n_k}\to\O_{\Sigma_X}\to 0
$$
by $a$ and taking cohomology, we get the exact sequence
$$
H^s(\O_{\Sigma_X}(a))\lra H^{s+1}(\I_{\Sigma_X}(a))\lra H^{s+1}(\O_{\pr^n_k}(a)).
$$
The first and last groups above are zero because $\dim(\Sigma_X)<s$ and $r>0$, 
respectively. Thus 
$H^{s+1}(\I_{\Sigma_X}(a))=0$. 

So the boundary map in (\ref{midmap}) is surjective. Furthermore, since (\ref{s+1a}) holds, we have 
that 
$H^s(\I_{\Sigma_X,X}(a))\cong k$ if and only if the boundary map is injective, which is the 
case if and only if 
$H^s(\I_{\Sigma_X}(a))=0$. But, if $\sigma=0$ then $\Sigma_X=\emptyset$, and hence 
$\I_{\Sigma_X}=\O_{\pr^n_k}$; then $H^s(\I_{\Sigma_X}(a))=0$ because $s<n$. And if 
$r-2\geq\sigma$, then $a\geq\sigma-s$, and thus 
$H^s(\I_{\Sigma_X}(a))=0$.
\edem

\begin{obs}\normalfont The above proof establishes an equivalence:
$$
H^s(\Omega^s_X)\cong k\quad\text{if and only if}\quad
H^s(\I_{\Sigma_X}(r-s-2))=0.
$$
If $s=1$ then $\Sigma_X$ is finite. If $X$ is a line 
then $\sigma=0$. Otherwise, $r\geq 2$, and it follows from 
Proposition~\ref{regconjfin} that $H^1(\I_{\Sigma_X}(r-3))=0$ only 
if $\sigma\leq r-2$. In other 
words, the converse to Theorem \ref{h11} holds if $s=1$.
\end{obs}

\section{Pfaff fields} 

\begin{subsct}\normalfont (\textit{Pfaff fields}) 
Let $V$ be an algebraic $k$-scheme. 
By definition, a \textit{Pfaff field} on $V$ is a map 
$\eta:\Omega_V^s\rightarrow\L$ of $\O_V$-modules, where $\L$ is an invertible sheaf on $V$ and 
$s$ is a positive integer. We call $s$ the \textit{rank} of $\eta$. Define the \textit{singular locus}
of $\eta$ to be the closed subscheme $\S\subseteq V$ defined by the sheaf of ideals
$\textrm{Im}(\eta\otimes\L^{-1})$. 

A closed subscheme $X\subseteq V$ is said to be
\textit{invariant} under $\eta$ if there is a Pfaff field 
$\varphi:\Omega_X^s\rightarrow\L|_X$ making the following
diagram commute:
$$
\begin{CD}
\Omega^s_V @>\eta >> \L\\
@VVV @VVV\\
\Omega^s_X @>\varphi >> \L|_X,
\end{CD}
$$
where the vertical maps are the natural restrictions.

If $X\subseteq V$ is reduced and invariant by $\eta$, then any union $Y$ of components of $X$, 
with its reduced induced subscheme structure, is also invariant by $\eta$. Indeed, in this situation, the 
restriction $\Omega^s_X|_Y\to\Omega^s_Y$ is surjective with generically zero kernel, 
and thus any map $\Omega^s_X|_Y\to\L|_Y$ factors through the restriction.

If $V=\pr^n_k$, and $\eta\:\Omega^s_{\pr^n_k}\to\L$ is a nonzero Pfaff field on 
$\pr^n_k$ of rank $s<n$, then $m\geq 0$, where $m:=\deg(\L)+s$. 
Indeed, since $\pr^n_k$ is smooth of 
dimension $n$, the field $\eta$ corresponds to a nonzero element of
$$
H^0(\Omega^{n-s}_{\pr^n_k}\otimes\L\otimes(\Omega^n_{\pr^n_k})^{-1}).
$$
So $H^0(\Omega^{n-s}_{\pr^n_k}(m+n+1-s))\neq 0$. 
By \cite{D}, Thm.~1.1, p.~40, 
this is only possible if $m+n+1-s>n-s$, that is, if $m\geq 0$. 
\end{subsct}

\begin{teo}
\label{principal0} Let $X\subseteq\pr^n_k$ be a connected, reduced, arithmetically Cohen--Macaulay 
subcanonical subscheme of pure dimension $s>0$ and degree $d$. Let $\Sigma_X$ be the singular 
locus of $X$. Assume the characteristic of $k$ is $0$ or does not divide $d$. 
Assume $X$ is invariant under a Pfaff field 
$\eta\:\Omega^s_{\pr^n_k}\to\L$ of rank $s$ in such a way that 
no irreducible component of $X$ is contained in the singular locus of $\eta$. Set 
$$
\sigma:=\reg(\Sigma_X),\quad r:=\reg(X),\quad m:=\deg(\L)+s.
$$
If $\sigma=0$ or $\sigma\leq r-2$ then $r\leq m+1$.
\end{teo}

\bdem By Theorem \ref{h11}, we have $h^s(\Omega^s_X)=1$. So, by \cite{EsK2}, Cor.~4.5, p.~3790, since 
$X$ is a.C.M., $r\leq s+\deg(\L)+1$, as claimed.
\edem

\begin{lema}\label{sigmaXY} Let $X$ be an equidimensional, reduced, projective $k$-scheme. Let 
$Y$ be a union of irreducible components of $X$, with its reduced induced subscheme structure. 
Let $\Sigma_X$ and $\Sigma_Y$ be the singular loci of $X$ and $Y$. Then $\Sigma_Y\subseteq\Sigma_X$.
\end{lema}

\begin{proof} 
Let $\mathcal H$ be the cokernel of $\gamma_X$ and $\mathcal G$ that of $\gamma_Y$. It is 
enough to observe that $\mathcal G$ is a subsheaf of a quotient of $\mathcal H|_Y$.  If $Y=X$ the 
assertion is trivial. So assume $Y\neq X$. Let $Z:=\overline{X-Y}$, again with the reduced induced 
subscheme structure. From the way $\gamma_X$ is defined, we see that 
$\gamma_X$ decomposes as
\begin{equation}\label{compog}
\begin{CD}
\Omega^s_X @>>> \Omega^s_Y\oplus\Omega^s_Z @>(\gamma_Y,\gamma_Z)>> \omega_Y\oplus\omega_Z @>\lambda>> 
\omega_X,
\end{CD}
\end{equation}
where the first map is induced by restriction of forms, and the last map, $\lambda$, is induced 
from the natural restriction map $\O_X\to\O_Y\oplus\O_Z$. Let $\T(\omega_X|_Y)$ be the torsion 
subsheaf of $\w_X|_Y$, and denote by $\w_{X,Y}$ the quotient. 
Restricting (\ref{compog}) to $Y$ and removing torsion, we get the 
following composition:
$$
\begin{CD}
\Omega^s_X|_Y @>\beta >> \Omega^s_Y @>\gamma_Y>> \omega_Y @>\iota >> \w_{X,Y},
\end{CD}
$$
where $\beta$ is the restriction map of $s$-forms, 
and $\iota$ is the composition of the canonical injection
$\w_Y\to\w_Y\oplus\w_Z$ with $\lambda$ and the quotient map $\w_X\to\w_{X,Y}$. Since 
$\lambda$ is generically an isomorphism and $\w_Y$ is torsion-free, $\iota$ is injective. Since 
$\beta$ is surjective and $\iota$ is injective, we get an injective map from $\mathcal G$ to 
$$
\frac{\omega_X|_Y}{\Imag(\gamma_X|_Y)+\T(\omega_X|_Y)},
$$
which is a quotient of $\mathcal H|_Y$.
\end{proof} 

\begin{teo}
\label{principal1} Let $X\subseteq\pr^n_k$ be a reduced, arithmetically Cohen--Macaulay 
subscheme of pure dimension $s>0$. Let $\Sigma_X$ be the singular locus of $X$. 
Assume $X$ is invariant under a Pfaff field $\eta\:\Omega^s_{\pr^n_k}\to\L$ of rank $s$ in such a way that 
no irreducible component of $X$ is contained in the singular locus of $\eta$. Set
$$
\sigma:=\reg(\Sigma_X),\quad r:=\reg(X),\quad m:=\deg(\L)+s.
$$
Then $r\leq m+1+\rho$, where $\rho:=\max(1,\sigma-r+2)$.
\end{teo}

\bdem Set $a:=r-s-2$. Let $\ell$ be any integer such that $\ell\geq\rho$. Since 
$\rho\geq 1$, we have $a+\ell\geq r-s-1$. Let $H\subset\pr^n_k$ be a general hyperplane. Multiplication 
by $(a+\ell-r+s+1)H$ induces an injection $\O_X(r-s-1)\to\O_X(a+\ell)$. Since 
$H^s(\O_X(r-s-1))=0$ by Proposition~\ref{reginter1}, and since the cokernel of the injection 
is supported in dimension at most $s-1$, it follows that 
\begin{equation}\label{H1Yrho}
H^s(\O_X(a+\ell))=0.
\end{equation}
Let $\I_X$ and $\I_{\Sigma_X}$ be the sheaves of ideals of $X$ and $\Sigma_X$ in $\pr^n_k$. 
Twisting the natural short exact sequence
$$
0\to\I_X\to\O_{\pr^n_k}\to\O_X\to0
$$
by $a+\ell$ and taking cohomology we get the exact sequence
$$
H^s(\O_X(a+\ell))\to H^{s+1}(\I_X(a+\ell))\to H^{s+1}(\O_{\pr^n_k}(a+\ell)).
$$
Since $a+\ell\geq r-s-1>-s-1\geq-n-1$, the last group is zero by Serre computation, and thus, 
using (\ref{H1Yrho}), we get 
\begin{equation}\label{H2Y0}
H^{s+1}(\I_X(a+\ell))=0.
\end{equation}
On the other hand, since $\sigma=\reg(\Sigma_X)$ and $a+\rho\geq\sigma-s$, we have
\begin{equation}\label{H2SY0}
H^s(\I_{\Sigma_X}(a+\ell))=0.
\end{equation}

If $Y$ is a union of irreducible components of $X$, with its reduced induced subscheme structure, 
then, since $\I_X\subset\I_Y$ with quotient supported in dimension at most $s$, Equation~\eqref{H2Y0} 
implies that
\begin{equation}\label{H2YY}
H^{s+1}(\I_Y(a+\ell))=0.
\end{equation}
Similarly, since $\I_{\Sigma_X}\subset\I_{\Sigma_Y}$ by Lemma \ref{sigmaXY}, and the quotient is 
supported in dimension at most $s-1$, Equation \eqref{H2SY0} implies
\begin{equation}\label{H2SYY0}
H^s(\I_{\Sigma_Y}(a+\ell))=0.
\end{equation}

Twisting the short exact sequence
\begin{equation}\label{exSingC}
0\to\I_Y\to\I_{\Sigma_Y}\to\I_{\Sigma_Y,Y}\to 0
\end{equation} 
by $a+\ell$, and taking cohomology, we get
the exact sequence 
$$
H^s(\I_{\Sigma_Y}(a+\ell))\to H^s(\I_{\Sigma_Y,Y}(a+\ell))\to H^{s+1}(\I_Y(a+\ell)).
$$
Using (\ref{H2YY}) and (\ref{H2SYY0}) we get that 
\begin{equation}\label{SXX}
H^s(\I_{\Sigma_Y,Y}(a+\ell))=0.
\end{equation}

Now, it follows from Proposition \ref{reginter1} that $H^s(\O_X(a))\neq 0$. Thus, by Serre Duality, 
there is a nonzero map $\tau\:\O_X(a)\to\w_X$. If $X$ is subcanonical, this map is an isomorphism. At any 
rate, since both $\O_X(a)$ and $\w_X$ are torsion-free, there is a union 
$Y$ of irreducible components of $X$, with its reduced induced subscheme structure, such that 
$\tau$ factors though an injection $\O_Y(a)\to\w_X$. This map factors through the natural map 
$\w_Y\to\w_X$, yielding an injection $\O_Y(a)\to\w_Y$. Of course, this injection induces one from 
$\I_{\Sigma_Y,Y}(a)$ to $\I_{\Sigma_Y,Y}\w_Y$, which can be composed with the injection 
$\I_{\Sigma_Y,Y}\w_Y\to\wt\Omega^s_Y$, where 
$$
\wt\Omega^s_Y:=\frac{\Omega^s_Y}{\T(\Omega^s_Y)}, 
$$
with $\T(\Omega^s_Y)$ denoting the torsion subsheaf of $\Omega^s_Y$. Since $\I_{\Sigma_Y,Y}(a)$ and 
$\wt\Omega^s_Y$ are rank-1, the cokernel of the composition $\I_{\Sigma_Y,Y}(a)\to\wt\Omega^s_Y$ 
is supported 
in dimension at most $s-1$. Thus, it follows from \eqref{SXX} that
\begin{equation}\label{SYY}
H^s(\wt\Omega^s_Y(\ell))=0.
\end{equation}

Notice that $\L=\O_{\pr^n_k}(m-s)$. 
Since $X$ is invariant under $\eta$, so is $Y$. So there is a Pfaff field 
$\varphi\:\Omega^s_Y\to\O_Y(m-s)$ making 
the following diagram commute:
$$
\begin{CD}
\Omega^s_{\pr^n_k} @>\eta >> \O_{\pr^n_k}(m-s)\\
@VVV @VVV\\
\Omega^s_Y @>\varphi >> \O_Y(m-s)
\end{CD}
$$
where the vertical maps are the natural restrictions. The image of $\eta$ is by definition 
$\I_{\S,\pr^n_k}(m-s)$, where $\S$ is the singular locus of $\eta$. So, since the 
vertical maps are surjective, the image of $\varphi$ is $\I_{\S\cap Y,Y}(m-s)$.

Now, since $\dim(\S\cap Y)<s$, the map $\varphi$ is generically surjective, and hence, since 
$Y$ is generically smooth, generically injective. In this case, the kernel of $\varphi$ is 
the torsion subsheaf $\T(\Omega_Y^1)$. So $\wt\Omega^s_Y\cong\I_{\S\cap Y,Y}(m-s)$, 
and hence \eqref{SYY} implies that 
\begin{equation}\label{H1ISY}
H^s(\I_{\S\cap Y,Y}(m-s+\ell))=0.
\end{equation}

Twisting the natural exact sequence
$$
0\to\I_{\S\cap Y,Y}\to\O_Y\to\O_{\S\cap Y}\to0
$$
by $m-s+\ell$ and taking cohomology, we get the exact sequence:
$$
H^s(\I_{\S\cap Y,Y}(m-s+\ell))\to H^s(\O_Y(m-s+\ell))\to H^s(\O_{\S\cap Y}(m-s+\ell)).
$$
Since $\dim(\S\cap Y)<s$, the last group is zero. So, it follows from (\ref{H1ISY}) 
that
\begin{equation}\label{ellrho}
H^s(\O_Y(m-s+\ell))=0.
\end{equation}
However, since there is an injection $\O_Y(a)\to\w_Y$, we have that $H^s(\O_Y(a))\neq 0$. Since 
\eqref{ellrho} holds for each $\ell\geq\rho$, we have $a\leq m-s+\rho-1$, 
from which follows the stated inequality.
\edem

\begin{cor}\label{regfolcur}
Let $X\subseteq\pr^n_k$ be a reduced complete intersection of hypersurfaces 
of degrees $d_1,\dots,d_{n-s}$ for a certain positive integer $s$. 
Let $\Sigma_X$ be the singular locus of $X$. Set $\sigma:=\reg(\Sigma_X)$ and put
$$
\rho:=\sigma+n-s+1-d_1-\cdots-d_{n-s}.
$$
Assume the characteristic of $k$ is $0$ or does not divide any of the $d_i$. 
Assume $X$ is invariant under a Pfaff field $\eta\:\Omega^s_{\pr^n_k}\to\L$ of rank $s$. Set 
$m:=\deg(\L)+s$. If $\dim(\S\cap X)<s$, where
$\S$ is the singular locus of $\eta$, then
$$
d_1+\dots+d_{n-s}\leq\left\{\begin{array}{lc}
 m+n-s,       & \mbox{ if } \rho\leq 0,  \\
 m+n-s+\rho,  & \mbox{ if } \rho > 0.
\end{array}\right.
$$
\end{cor}

\bdem Since $X$ is a complete intersection, and of positive dimension, $X$ is a.C.M.~and connected. 
Also, the conormal sheaf $\mathcal C$ of $X$ satisfies
$$
\mathcal C\cong\O_X(d_1)\oplus\cdots\oplus\O_X(d_{n-s}).
$$
Thus
$$
\omega_X\cong\bigwedge^n\Omega^1_{\pr^n_k}|_X\otimes(\bigwedge^{n-s}\mathcal C)^{\vee}\cong
\O_X(d_1+\cdots+d_{n-s}-n-1).
$$
Hence, by Proposition \ref{regsub},
$$
r=d_1+\cdots+d_{n-s}-n+s+1.
$$
Apply Theorems \ref{principal0} and \ref{principal1} now.
\edem

\section{Rank-1 Pfaff fields} 

\begin{subsct}\normalfont (\textit{Degree}) 
Let $\eta\:\Omega^1_{\pr^n_k}\to\L$ be a rank-1 Pfaff field. Set $m:=\deg(\L)+1$; then 
$\L\cong\O_{\pr^n_k}(m-1)$. 
By convention, we say that $\eta$  has \textit{degree} $m$. More geometrically, the degree of $\eta$ is 
the 
degree of the subscheme of points of a general hyperplane $H$ where the direction given by $\eta$ is 
contained in $H$. More precisely, given $H$, the subscheme is the degeneration scheme of the map of 
locally free sheaves
$$
\begin{CD}
\Omega^1_{\pr^n_k}|_H @>(\eta|_H,\beta)>> \O_H(m-1)\oplus\Omega^1_H,
\end{CD}
$$
where $\beta$ is the natural restriction. That the degree of this subscheme is indeed $m$ follows by 
taking determinants, noticing that $\det\Omega^1_{\pr^n_k}\cong\O_{\pr^n_k}(-n-1)$ and 
$\det\Omega^1_H\cong\O_H(-n)$.
\end{subsct} 

\begin{prop}\label{regsing} 
Let $\eta\:\Omega_{\pr ^n_k}^1\to\O_{\pr^n_k}(m-1)$ be a rank-$1$ Pfaff field on 
$\pr^n_k$, for $n\geq 2$. If $m\geq 1$ and the singular locus $\S$ of $\eta$ is finite then 
$$
\reg(\S)=nm-n+2.
$$
\end{prop}

\bdem Let $\I_\S$ be the sheaf of ideals of $\S$ and 
$\eta':=\eta(1-m)$. Then $\eta'\:\Omega^1_{\pr^n_k}(1-m)\to\O_{\pr^n_k}$ has image $\I_\S$, or 
degeneration scheme $\S$. Consider the Koszul complex of $\eta'$:
\begin{equation}\label{1}
\begin{CD}
0\to\Omega_{\pr ^n_k}^n(n-nm) @>d_n>> \cdots @>d_2>> \Omega_{\pr ^n_k}^1(1-m) @>d_1>> \O_{\pr^n_k}\to 0,
\end{CD}
\end{equation}
where $d_1:=\eta'$.
Since $\S$ is finite, $\S$ is of the expected codimension. Since $\pr^n_k$ is Cohen--Macaulay, 
the dual to $\eta'$ is a regular section, and hence the complex above is exact at positive level. 

Let $\I_j:=\Imag(d_j)$ for $j=1,\ldots,n$. Then $\I_1=\I_\S$ and 
$\I_n\cong\Omega_{\pr ^n_k}^n(n-nm)$. Also, we can break (\ref{1}) in the following short exact sequences:
\begin{equation} \label{sqrecorrencia}
 0\rightarrow
\I_{j+1} \rightarrow\Omega_{\pr ^n_k}^j(j-j m) \rightarrow
\I_j\rightarrow 0, \,\,\,\,\,\,\, j=1,\ldots,n-1.
\end{equation}
Twisting these sequences by $r-1$, and taking cohomology, we get the exact sequences:
\begin{equation}\label{exatas1}
H^j(\Omega_{\pr^n_k}^{j}(\ell_{j,r}))\lra H^j(\I_j(r-1))\lra H^{j+1}(\I_{j+1}(r-1))\lra 
H^{j+1}(\Omega_{\pr^n_k}^{j}(\ell_{j,r}))
\end{equation}
for $j=1,\ldots, n-1$, where $\ell_{j,r}:=j-jm+r-1$. 

Set $b:=nm-n+2$. Notice that, since $m\geq 1$,
$$
r-m=\ell_{1,r}\geq\ell_{2,r}\geq\cdots\geq\ell_{n-1,r}=r+m-b.
$$
So, if $r\geq b-1$ then $\ell_{j,r}\geq m-1$ for $j=1,\dots,n-1$, with equality only if 
$r=b-1$. In particular, $\ell_{j,r}\geq 0$ for $r\geq b-1$, with equality only if 
$r=b-1$. So, from \cite{D}, Thm.~1.1, p.~40, it follows that 
$H^j(\Omega_{\pr^n_k}^{j}(\ell_{j,r}))=0$ for $r\geq b$, while 
$H^{j+1}(\Omega_{\pr^n_k}^{j}(\ell_{j,r}))=0$ for $r\geq b-1$, for $j=1,\ldots, n-1$. 
Then, from the exact sequences (\ref{exatas1}) we get surjections
$$
H^1(\I_1(r-1))\lra H^{2}(\I_2(r-1))\lra\cdots\lra H^{n-1}(\I_{n-1}(r-1))\lra H^n(\I_n(r-1))
$$
for $r\geq b-1$, which are all isomorphisms for $r\geq b$. 
Now, $\I_n(r-1)\cong\Omega^n_{\pr^n_k}(n-nm+r-1)$. So, again by \cite{D}, Thm.~1.1, p.~40, we have that
$h^n(\I_n(r-1))\neq 0$ if $r\leq b-1$, whereas $h^n(\I_n(r-1))=0$ if $r\geq b$. 
Then $\reg(\S)=b$ by Proposition \ref{regconjfin}.
\edem

\begin{obs}\normalfont If $m=0$ and $\eta\neq 0$ then $\S$ consists of a point, and thus $\reg(\S)=1$.
\end{obs}

\begin{teo}\label{igualdade1} Let $C\subseteq\pr^n_k$ be a reduced, arithmetically 
Cohen--Macaulay, subcanonical subscheme of dimension $1$. Let $\Sigma_C$ be the singular locus of $C$. 
Assume $C$ is invariant under a rank-$1$ Pfaff field $\eta\:\Omega^1_{\pr^n_k}\to\L$ of degree $m\geq 1$. 
Set
$$
\sigma:=\reg(\Sigma_C)\quad\text{and}\quad r:=\reg(C).
$$
Assume that $r\geq 5$ if $m=1$ or $r\geq mn-n+4$ if $m>1$. If the singular locus of $\eta$ is finite , 
then $r=m+1+\rho$, where $\rho:=\sigma-r+2$.
\end{teo}

\bdem Since $r\geq 4$, we have $n\geq 2$. Then $r\geq m+4$. Indeed, 
$$
n(m-1)+4\geq 2(m-1)+4=m+(m-2)+4\geq m+4
$$
if $m>1$. So $\rho\geq 3$ and $r\leq m+1+\rho$ by Theorem \ref{principal1}. In particular, 
$\sigma>0$. We need only prove that $r\geq m+1+\rho$.

Let $\S$ denote the singular locus of $\eta$. 
Let $\I_\S$ and $\I_C$ be the sheaves of ideals of $\S$ and $C$, and $\I_{\S\cap C}$ that of 
$\S\cap C$ in $\pr^n_k$. Set $j:=m+\rho-2$. Twisting the natural short exact sequence
$$
0\to\I_\S\to\I_{\S\cap C}\to\I_{\S\cap C,\S}\to 0
$$
by $j$, and taking cohomology, we obtain the exact sequence
\begin{equation}\label{cohSC} 
H^1(\I_\S(j))\lra H^1(\I_{\S\cap C}(j))\lra H^1(\I_{\S\cap C,\S}(j)).
\end{equation}
Since $\S$ is finite, the last group is zero. Furthermore, since $r\leq m+1+\rho$, we have
$$
j+1\geq r-2\geq mn-n+2.
$$ 
Since $\reg(\S)=mn-n+2$ by Proposition \ref{regsing}, also $H^1(\I_\S(j))=0$. Thus
$$
H^1(\I_{\S\cap C}(j))=0.
$$

Now, twist the natural short exact
sequence
$$
0\to\I_C\to\I_{\S\cap C}\to\I_{\S\cap C,C}\to 0
$$
by $j$, and take cohomology to get the exact sequence
\begin{equation}\label{eql1} 
H^1(\I_{\S\cap C}(j))\lra H^1(\I_{\S\cap C,C}(j))\lra H^2(\I_C(j)).
\end{equation}
Since $H^1(\I_{\S\cap C}(j))=0$, if we show that $H^1(\I_{\S\cap C,C}(j))\neq 0$, then it 
follows from the exactness of \eqref{eql1} that $H^2(\I_C(j))\neq 0$, and hence that 
$r\geq j+3$. 

Since $j+3=m+\rho+1$, we need only show that $H^1(\I_{\S\cap C,C}(j))\neq 0$. Since 
$C$ is invariant under $\eta$, and $\S$ is finite, we have that 
$\I_{\S\cap C,C}(m-1)\cong\wt\Omega^1_C$, where 
$$
\wt\Omega^1_C:=\frac{\Omega^1_C}{\T(\Omega^1_C)},
$$
with $\T(\Omega^1_C)$ denoting the torsion subsheaf of $\Omega^1_C$. Since $C$ is subcanonical, 
$\w_C\cong\O_C(r-3)$ by Proposition \ref{regsub}. Furthermore, $C$ is Gorenstein, whence 
$\I_{\Sigma_C,C}\,\w_C\cong\wt\Omega^1_C$. So, since $j=m+\rho-2$ and $r+\rho-2=\sigma$, it follows that  
$H^1(\I_{\S\cap C,C}(j))\neq 0$ is equivalent to 
\begin{equation}\label{final}
H^1(\I_{\Sigma_C,C}(\sigma-2))\neq 0.
\end{equation}

Let $\I_{\Sigma_C}$ be the sheaf of ideals of $\Sigma_C$ in $\pr^n_k$. 
Twisting the natural exact sequence
$$
0\to\I_C\to\I_{\Sigma_C}\to\I_{\Sigma_C,C}\to 0
$$
by $\sigma-2$, and taking cohomology, we get the exact sequence
\begin{equation}\label{cohCsigma} 
H^1(\I_C(\sigma-2))\longrightarrow H^1(\I_{\Sigma_C}(\sigma-2))\longrightarrow 
H^1(\I_{\Sigma_C,C}(\sigma-2)).
\end{equation}
Since $r\geq m+4$, we have
$$
r\leq m+\rho+1=m+3+\sigma-r\leq m+3+\sigma-m-4=\sigma-1.
$$
So, since $r=\reg(C)$, we have
$$
H^1(\I_C(\sigma-2))=0.
$$
On the other hand, since $\Sigma_C$ is finite and nonempty, $H^1(\I_{\Sigma_C}(\sigma-2))\neq 0$ by 
Proposition~\ref{regconjfin}. So, from the exactness of (\ref{cohCsigma}) we get 
\eqref{final}.
\edem


\begin{thebibliography}{CmCrG}

\bibitem[BMe]{BMe} M. Brunella and L. G. Mendes,
\emph{Bounding the degree of solutions to Pfaff equations}, 
Publ. Mat. {\bfseries 44} (2000), 593--604.

\bibitem[CmCr]{CamCar} A. Campillo  and  M. Carnicer,
\emph{Proximity inequalities and bounds for the degree of invariant curves by foliations of 
$\pr^2_{\complexo}$}.
Trans. Amer. Math. Soc. {\bfseries 349} (1997), 2211--2228.

\bibitem[CmCrG]{CamCarG} A. Campillo  and  M. Carnicer and J. Garc\'{i}a de la Fuente,
\emph{Invariant curves by vector fields on algebraic varieties}.
J. London Math. Soc. (2) {\bfseries 62} (2000), 56--70.

\bibitem[Cr]{Car} M. Carnicer,
\emph{The Poincar\'e problem in the nondicritical case}.
Ann. Math. {\bfseries 140} (1994) 289--294.

\bibitem[CeLn]{CeL} D. Cerveau and A. Lins Neto,
\emph{Holomorphic foliations in $CP(2)$ having an invariant algebraic curve}.
Ann. Inst. Fourier {\bfseries 41} (1991), 883--903.

\bibitem[D]{D} P. Deligne,
\emph{Cohomologie des intersections completes} (SGA 7 II).
Lecture Notes in Mathematics, vol. 340, pp. 39--61, Springer-Verlag, 1973.

\bibitem[E]{Es} E. Esteves,
\emph{The Castelnuovo--Mumford regularity of an integral variety of a vector field on projective space}.
Math. Res. Letters {\bfseries 9} (2002), 1--15.

\bibitem[EKl1]{EsK1} E. Esteves and S. Kleiman,
\emph{Bounds on leaves of one-dimensional foliations}.
Bull. Braz. Math. Soc. {\bfseries 34} (2003), 145--169.

\bibitem[EKl2]{EsK2} E. Esteves and S. Kleiman,
\emph{Bounding solutions of Pfaff equations}.
Comm. Algebra {\bfseries 31} (2003), 3771--3793.

\bibitem[EKl3]{EsK3} E. Esteves and S. Kleiman,
\emph{Bounds on leaves of foliations of the plane}.
Contemp. Math. {\bfseries 354} (2004), 57--67.

\bibitem[J]{J} J. P. Jouanolou,
\emph{Equations de Pfaff alg\'ebriques}.
Lecture Notes in Mathematics, vol. 708, Springer-Verlag, 1979.

\bibitem[Ku]{Ku} E. Kunz,
\emph{Holomorphe Differentialformen auf algebraischen Variet\"aten 
mit Singularit\"aten}.
Manuscripta math. {\bfseries 25} (1975), 91--108.

\bibitem[Ln]{L} A. Lins Neto,
\emph{Some examples for Poincar\'e and Painlev\'e problems}. 
Ann. Scient. \'Ec. Norm. Sup. {\bfseries 35} (2002), 231--266.

\bibitem[Lp]{Lp} J. Lipman,
\emph{Dualizing sheaves, differentials and residues on
algebraic varieties}.
Ast\'erisque, vol.~117, Soc. Math. de France, 1984.

\bibitem[dPW]{dPW} A. A. du Plessis and C. T. C. Wall,
\emph{Application of the theory of the discriminant to highly singular plane curves}.
Math. Proc. Camb. Phill. Soc. {\bfseries 126} (1999), 259--266.

\bibitem[Mu]{Mu} D. Mumford,
\emph{Lectures on curves on an algebraic surface}.
Annals of Mathematical Studies, vol. 59, Princeton University Press, Princeton, 
1966.

\bibitem[Pe]{Pe} J. V. Pereira,
\emph{On the Poincar\'e problem for foliations of general type}.
Math. Ann. {\bfseries 323} (2002), 217--226.

\bibitem[Po]{Po} H. Poincar\'e,
\emph{Sur l'int\'egration alg\'ebrique des \'equations differentielles du premier ordre et du premier 
degr\'e}.
Rendiconti del Circolo Matematico di Palermo {\bfseries 5} (1891), 161--191.

\bibitem[S]{S} M. Soares, 
\emph{The Poincar\'e problem for hypersurfaces invariant by one-dimensional 
foliations}. 
Invent. Math. {\bfseries 128} (1997), 495--500.

\end{thebibliography}
\end{document}